\documentclass[12pt]{article}
\usepackage{amsmath}
\usepackage{amssymb}
\usepackage{amscd}
\usepackage{latexsym}
\usepackage{theorem}
\usepackage{doublespace}
\usepackage{epsfig}
\usepackage{psfrag}
\usepackage{amsfonts}
\usepackage{enumerate}

\newtheorem{theorem}{Theorem}[section]
\newtheorem{lemma}[theorem]{Lemma}
\newtheorem{proposition}[theorem]{Proposition}
\newtheorem{corollary}[theorem]{Corollary}
\newtheorem{claim}[theorem]{Claim}

{\theorembodyfont{\rmfamily}
\theoremstyle{plain}

\newtheorem{example}[theorem]{Example}

\newtheorem{remark}[theorem]{Remark}
}

\newcommand{\HH}{{\mathbb H}}

\newcommand{\liet}{{\mathfrak t}}

\newcommand{\e}{\varepsilon}

\renewcommand{\k}{\kappa}
\newcommand{\ktd}{\kappa_{\Td}}
\newcommand{\ktdso}{\kappa_{\Td\times S^1}}
\newcommand{\kso}{\kappa_{S^1}}

\newcommand{\<}{\left<}
\renewcommand{\>}{\right>}

\renewcommand{\dim}{\operatorname{dim}}

\renewcommand{\H}{{\mathbb{H}}}
\newcommand{\C}{{\mathbb{C}}}
\newcommand{\Zt}{{\mathbb{Z}_2}}
\newcommand{\Z}{{\mathbb{Z}}}
\newcommand{\Q}{{\mathbb{Q}}}
\newcommand{\R}{{\mathbb{R}}}

\renewcommand{\a}{\alpha}

\newcommand{\otd}{\{1,\ldots,d\}}
\newcommand{\Hn}{{\H^n}}
\newcommand{\Cn}{{\C^n}}

\renewcommand{\cot}{T^*\Cn}
\newcommand{\Tk}{T^k}
\newcommand{\Tn}{T^n}

\newcommand{\Tdr}{T_{\R}^d}
\newcommand{\Td}{T^d}

\newcommand{\Tnd}{(\tn_{\Z})^*}

\newcommand{\tk}{\mathfrak{t}^k}
\newcommand{\tn}{\mathfrak{t}^n}
\newcommand{\td}{\mathfrak{t}^d}
\newcommand{\tkd}{(\tk)^*}
\newcommand{\tnd}{(\tn)^*}
\newcommand{\tdd}{(\td)^*}
\newcommand{\Tdso}{\Td\times S^1}

\newcommand{\subs}{\subseteq}

\newcommand{\hookto}{{\hookrightarrow}}
\newcommand{\onto}{{\twoheadrightarrow}}
\renewcommand{\iff}{\Leftrightarrow}
\newcommand{\impl}{\Rightarrow}

\newcommand{\hm}{H^*(M)}
\newcommand{\htm}{H^*_{\Td}(M)}
\newcommand{\htsm}{H^*_{\Td\times S^1}(M)}
\newcommand{\hsm}{H^*_{S^1}(M)}
\newcommand{\hr}{H^*_{\Tdr}(M_{\R};\Zt)}
\newcommand{\hrs}{H^*_{\Tdr\times\Zt}(M_{\R};\Zt)}
\newcommand{\hscomp}{H^*_{\Zt}(\comp;\Zt)}
\newcommand{\hcomp}{H^*(\comp;\Zt)}

\newcommand{\hs}{\hspace{3pt}}

\newcommand{\gd}{\g^*}
\newcommand{\g}{\mathfrak{g}}
\newcommand{\half}{\frac{1}{2}}
\renewcommand{\mod}{{/\!\!/}}
\newcommand{\mmod}{{/\!\!/\!\!/\!\!/}}
\newcommand{\mc}{{\mathcal{C}}}
\newcommand{\arr}{{\mathcal{H}}}
\newcommand{\comp}{\mathcal{M}(\arr)}
\newcommand{\bomp}{\Cn\setminus\cup_{i=1}^n H_i^{\C}}
\newcommand{\os}{\mathcal{OS}}
\newcommand{\ost}{\mathcal{OS}\otimes\Zt}
\newcommand{\mhk}{\mu_{\text{HK}}}
\renewcommand{\mc}{\mu_{\C}}
\newcommand{\mr}{\mu_{\R}}

\newcommand{\qed}{\hfill \mbox{$\Box$}\medskip\newline}
\newenvironment{proofc}{\noindent {\bf Proof of \ref{coorientation}:}}{\qed\par}
\newenvironment{proofp}{\text{ }\newline\noindent 
{\bf Proof of \ref{placement}:}}{\qed \par}
\newenvironment{proofhtsm}{\text{ }\newline\noindent 
{\bf Proof of \ref{htsm}:}}{\qed \par}
\newenvironment{proofhsm}{\text{ }\newline\noindent 
{\bf Proof of \ref{hsm}:}}{\qed \par}
\newenvironment{proof}{\noindent {\bf Proof:}}{\qed \par}

\newenvironment{sketch}{\noindent {\bf Sketch of Proof:}}{\qed \par}

\begin{document}
\begin{spacing}{1.1}

\noindent
{\LARGE \bf Properties of 
the residual circle action \smallskip \\
on a hypertoric variety}\bigskip\\
{\bf Megumi Harada} \\
Department of Mathematics, University of California,
Berkeley, CA 94720\smallskip \\
{\bf Nicholas Proudfoot } \\
Department of Mathematics, University of California,
Berkeley, CA 94720
\bigskip
{\small
\begin{quote}
\noindent {\em Abstract.}
We 
consider an orbifold $X$ obtained by a 
K\"ahler reduction of $\Cn$,
and we define its ``hyperk\"ahler analogue'' $M$
as a hyperk\"ahler reduction of $\cot\cong\Hn$ by the same group.
In the case where the group is abelian and $X$ is a toric variety,
$M$ is a toric hyperk\"ahler orbifold, as defined in
\cite{BD}, and further studied in \cite{K1,K2} and \cite{HS}.  
The variety $M$ carries a natural action of $S^1$, induced
by the scalar action of $S^1$ on the fibers of $\cot$.
In this paper we study this action, computing its fixed points
and its equivariant cohomology.
As an application, we use the associated $\Zt$ action
on the real locus of $M$ to compute a deformation of the
Orlik-Solomon algebra of a smooth, real hyperplane arrangement $\arr$,
depending nontrivially on the affine structure of the arrangement.
This deformation is given
by the $\Zt$-equivariant cohomology of the complement of the complexification
of $\arr$, where $\Zt$ acts by complex conjugation.
\end{quote}
}
\bigskip

In order to construct a toric variety as a K\"ahler quotient of $\Cn$ by a 
torus,
one begins with the combinatorial data of an arrangement $\arr$ of $n$ cooriented,
rational, affine hyperplanes in $\R^d$.
The normal vectors to these hyperplanes determine
a subtorus $\Tk\subs\Tn$ ($k=n-d$), and the affine structure determines a 
value $\a\in\tkd$ at which to reduce, so that we may define $X=\Cn\mod_{\!\!\a}\Tk$.
Using the same combinatorial data, one can also construct a {\em hypertoric 
variety},\footnote{In 
\cite{BD,K1,K2,HS} $M$ is called a ``toric hyperk\"ahler'' variety,
but as it is a complex variety that is not toric in
the standard sense, we prefer the term ``hypertoric.''}
which is defined as the hyperk\"ahler quotient $M = \Hn\mmod_{\!\!(\a,0)}\Tk$ of 
$\Hn\cong\cot$ by the
induced action of the same subtorus $\Tk\subs\Tn$ \cite{BD}.  It is well known 
that the toric variety $X$ does not retain all of the information of $\arr$; indeed, it 
depends only on the polyhedron
$\Delta$ obtained by intersecting the half-spaces associated to each of the 
cooriented hyperplanes. Thus it is always possible to add an extra hyperplane to $\arr$ 
without changing $X$.  In contrast, the hypertoric variety $M$ remembers the 
number of hyperplanes in $\arr$, but its equivariant diffeomorphism type 
depends neither on the coorientations nor on the affine structure of $\arr$ (see Theorem \ref{htm}, 
and Lemmas \ref{placement} and \ref{coorientation}).

The purpose of this paper is to study the hamiltonian $S^1$ action on $M$ descending
from the scalar action of $S^1$ on the fibers of $\cot$.  This action is sensitive to
both the coorientations and the affine structure of $\arr$,
even on the level of equivariant cohomology (Section \ref{Konno}).
One can recover the toric variety $X$ as the minimum of the $S^1$ moment map,
hence the geometric structure of $M$ along with its circle action carries strictly more
information than either $X$ or $M$ alone.
In Section \ref{coreflow} we give an explicit description of this action
when restricted to the core $C$, a deformation retract of
$M$ which is a union of projective subvarieties.
In Section \ref{Konno} we compute the $S^1$ and $\Tdso$-equivariant cohomologies
of $M$, using the full combinatorial data of $\arr$. 

In Section \ref{defos}, we 
examine the {\em real locus} $M_{\R}\subs M$, i.e. the fixed point set 
of an involution
of $M$ that is anti-holomorphic with respect to the first complex structure.
By studying the topology of $M_{\R}$, we interpret the results of Section \ref{Konno}
in terms of the Orlik-Solomon algebra 
$\mathcal{OS}=H^*(\comp)$,
where $\comp$
is the complement of the complexification of $\arr$.
We show how to interpret Theorem \ref{htsm} as a computation
of $H^*_{\Zt}(\comp;\Zt)$,
a deformation of the Orlik-Solomon algebra of a smooth, real arrangement
that depends nontrivially on the affine structure.\footnote{A more
general computation of $H^*_{\Zt}(\comp;\Zt)$, in which $\arr$
is not assumed to be simple, rational, or smooth, will appear in \cite{Pr}.}
\newline

\noindent {\em Acknowledgments.}
We are very grateful to Tam\'as Hausel for introducing us to this problem
and answering many questions.
We would also like to thank Allen Knutson for guiding our research.

\begin{section}{Hyperk\"ahler reductions}\label{reduction}
A hyperk\"ahler manifold is a smooth manifold, 
necessarily of real dimension $4n$,
which admits
three complex structures $J_1,J_2,J_3$
satisfying the usual quaternionic relations, in a manner compatible
with a metric. 
Just as in the K\"ahler case, we can define three different
symplectic forms on $N$ as follows:
\[
\omega_1(v,w) = g(J_{1}v, w), \hs\omega_2(v,w) = g(J_{2}v,w), \hs\omega_3(v,w) =
g(J_{3}v,w).
\]
Note that the complex-valued two-form \(\omega_2 + i
\omega_3\) is non-degenerate and covariant constant, hence closed and
{\em holomorphic} with respect to the complex structure $J_1$. 
Any hyperk\"aher manifold can therefore be considered as a holomorphic
symplectic manifold with complex structure $J_1$, real symplectic
form $\omega_{\R} := \omega_1$, and holomorphic symplectic form
$\omega_{\C} := \omega_2 + i\omega_3$.
This is the point of view that we will adopt in this paper.

We will refer to an action of $G$ on a hyperk\"ahler manifold
$N$ as {\em hyperhamiltonian}
if it is hamiltonian with respect to $\omega_{\R}$ and holomorphic
hamiltonian with respect to $\omega_{\C}$, with $G$-equivariant
moment map
$$\mhk:=\mr\oplus\mc:N\to\gd\oplus\g_{\C}^*.$$ 

\begin{theorem}{\em\cite{HKLR}}
Let \((N^{4n},g)\) be a hyperk\"ahler manifold with real symplectic form
$\omega_{\R}$ and holomorphic symplectic form $\omega_{\C}$.
Suppose that $N$ is equipped with a hyperhamiltonian action of a 
compact Lie group $G$, with moment map
$\mhk=\mr\oplus\mc$. 
Suppose $\xi = \xi_{\R}\oplus \xi_{\C}$ is a
central regular value of $\mhk$. Then there is a unique 
hyperk\"ahler structure
on the hyperk\"ahler quotient 
\(M = N \mmod_{\!\!\xi}G := \mhk^{-1}(\xi)/G\), 
with associated symplectic and holomorphic symplectic forms
$\omega^{\xi}_{\R}$ and $\omega^{\xi}_{\C}$, such that
$\omega^{\xi}_{\R}$ and $\omega^{\xi}_{\C}$ pull back to the restrictions
of $\omega_{\R}$ and $\omega_{\C}$ to $\mhk^{-1}(\xi)$.
\end{theorem}



If $\xi\in \gd\oplus\g_{\C}^*$ is fixed by the coadjoint action of $G$, the inverse
image $\mc^{-1}(\xi_{\C})$ is preserved by $G$,
and is a (singular) K\"ahler subvariety with respect to $\omega_{\R}$.
Then by 
\cite{HL} 
(see also \cite[3.2]{Na} and \cite[2.5]{Sj}),
we have 
$$N\mmod_{\!\!\xi}G = 
\mc^{-1}(\xi_{\C})\mod_{\!\!\xi_{\R}}G
=\mu_{\C}^{-1}(\xi_{\C})^{\text{ss}}/G_{\C},$$
where 
$$\mc^{-1}(\xi_{\C})^{\text{ss}} = \{x\in\mc^{-1}(\xi_{\C})\mid Gx\cap\mu_{\R}^{-1}(\xi_{\R})
\neq\emptyset\}.$$

We now specialize to the case where $G$ is a compact Lie group
acting linearly on $\Cn$ with moment map $\mu:\Cn \to \gd$, taking $0\in\C^n$
to $0\in\gd$.  This action 
induces an action of $G$ on the holomorphic cotangent bundle 
$\cot\cong\Cn\times(\Cn)^*$.
If we choose a bilinear inner product on $\Cn$, 
we can coordinatize this representation
as $\{(z,w)\mid z,w\in\Cn\}$ with $g (z,w)=(g z,g^{-1} w)$.
Choose an identification of $\Hn$ with $\cot$
such that the complex structure $J_1$ on $\Hn$ 
given by right multiplication by $i$
corresponds to the natural complex structure on $\cot$.
Then $\cot$ inherits a hyperk\"ahler, and 
therefore also a holomorphic symplectic, structure, 
with
$\omega_{\R}$ given by 
adding the standard symplectic structures on $\Cn$ and $(\Cn)^{*}
\cong \Cn$, and 
$\omega_{\C} = d\eta$, where $\eta$ is the canonical 
holomorphic 1-form on $\cot$.  


Note that $G$ acts $\H$-linearly on $\cot\cong\Hn$
(where $n\times n$ matrices
act on the left on $\Hn$, and scalar multiplication by $\H$
is on the right), 
and does so hyperhamiltonianly
with moment map $\mhk=\mc\oplus\mr$,
where
$$\mr(z,w) = \mu(z)-\mu(w)\hspace{10pt}\text{and}
\hspace{10pt}\mc(z,w)(v) = w(\hat{v}_z)$$
for $w\in T^*_z\Cn,\hs v\in\g_{\C}$, and $\hat{v}_z$ the element
of $T_z\Cn$ induced by $v$.
Consider a central regular value $\a\in\gd$ for $\mu$,
and suppose that $(\a,0)\in\gd\oplus\gd_{\C}$ is a central regular value for $\mu_{\text{HK}}$.
We refer to the hyperk\"ahler reduction $M = \Hn \mmod_{\!\!(\a,0)} G$
as the {\em hyperk\"ahler analogue} of the corresponding K\"ahler reduction
$X=\Cn\mod_{\!\!\a} G$. 
The following proposition is proven for the case 
where $G$ is a torus in \cite[7.1]{BD}. 

\begin{proposition}\label{compactification}
The cotangent bundle $T^*X$ is isomorphic to an open subset of $M$.
\end{proposition}

\begin{proof}
Let $Y = \{(z,w)\in\mc^{-1}(0)^{\text{ss}}\mid
z \in (\Cn)^{\text{ss}}\}$,
where we ask $z$ to be semistable with respect to $\a$
for the action of
$G_{\C}$ on $\Cn$, so that $X\cong (\Cn)^{\text{ss}}/G_{\C}$.
Let $[z]$ denote the element of $X$ represented by $z$.
The tangent space $T_{[z]}X$ is equal to the quotient of $T_{z}\Cn$
by the tangent space to the $G_{\C}$ orbit through $z$,
hence
$$T^*_{[z]}X\cong\{w\in T^*_{[z]}\Cn\mid w(\hat{v}_z)=0\text{ for all }
v\in\g_{\C}\} = \{w\in(\Cn)^*\mid\mc(z,w)=0\}.$$
Then $$T^*X \cong \{(z,w)\mid z\in(\Cn)^{\text{ss}}\text{ and }
\mc(z,w)=0\}/G_{\C} = Y/G_{\C}$$
is an open subset of $M$.
\end{proof}

Consider the action
of $S^1$ on $\Hn\cong T^*\Cn$ given by 
``rotating the fibers'' of the cotangent
bundle,  
given explicitly by $\tau(z,w) = (z,\tau w)$. 
This action is hamiltonian with respect to the real symplectic
structure $\omega_{\R}$ with moment map $\Phi(z,w) = \half |w|^2$.
Because it commutes with the action of $G$,
the action descends to a hamiltonian action on $M$, where we will
still denote the moment map by $\Phi$.
Since $S^1$ acts trivially on $z$, and by scalars on $w$,
it does {\em not} preserve the complex symplectic form
$\omega_{\C}(z,w) = dw\wedge dz$,
and does not act $\H$-linearly.

\begin{proposition}\label{proper}
If the original moment map $\mu:\Cn\to\gd$ is proper,
then so is $\Phi:M\to\R$.
\end{proposition}

\begin{proof}
We would like to show that $\Phi^{-1}[0,R]$ is compact
for any $R$.  Since 
$$\Phi^{-1}[0,R] = \{(z,w)\mid \mr(z,w)=\a,\hs
\mc(z,w)=0,\hs\Phi(z,w)\leq R\}\big/G$$ and $G$ is compact,
it is sufficient to show that 
the set $\{(z,w)\mid \mr(z,w)=\a,\hs\Phi(z,w)\leq R\}$
is compact.  Since $\mr(z,w)=\mu(z)-\mu(w)$,
this set is a closed subset
of $$\mu^{-1}\left\{\a+\mu(w) \,\big|\, \half |w|^2\leq R\right\}
\times\left\{w \,\big|\,
\half |w|^2 \leq R\right\},$$
which is compact by the properness of $\mu$.
\end{proof}


In the case where $G$ is abelian and $X$ is a nonempty toric 
variety, properness of $\mu$ (and therefore of $\Phi$) is equivalent 
to compactness of $X$. 
\end{section}

\begin{section}{Hypertoric varieties}\label{ht}
In this section we restrict our attention to 
{\em hypertoric}
varieties,
which are the hyperk\"ahler analogues
of toric varieties in the sense of
Section \ref{reduction}.  
We begin with the full $n$-dimensional torus $\Tn$ acting on $\Cn$,
and the induced action on $\Hn\cong\cot$
given by $t(z,w) = (tz, t^{-1} w)$.
Let \(\{a_i\}_{1 \leq i \leq n}\) be nonzero primitive integer vectors in
\(\liet^d \cong \R^d\) defining a map \(\beta: \liet^n \longrightarrow
\liet^d\) by \(\e_i \mapsto a_i,\) where
$\{\e_i\}$ is the standard basis for \(\liet^n \cong \R^n\), dual to $\{u_i\}$.
This map fits into an exact sequence
\[
0 \longrightarrow \liet^k \stackrel{\iota}{\longrightarrow} \liet^n
\stackrel{\beta}{\longrightarrow} \liet^d\longrightarrow 0,
\]
where \(\liet^k := {\mathrm ker}(\beta).\) Exponentiating, we
get the exact sequence
\[
0 \longrightarrow T^k \stackrel{\iota}{\longrightarrow} T^n
\stackrel{\beta}{\longrightarrow} T^d \longrightarrow 0,
\]
whereas by dualizing, we get
\[
0 \longrightarrow (\liet^d)^{*} \stackrel{\beta^{*}}{\longrightarrow}
(\liet^n)^{*} \stackrel{\iota^{*}}{\longrightarrow} (\liet^k)^{*} \longrightarrow 0,
\]
where we abuse notation by using $\iota$ and $\beta$
to denote maps on the level of groups as well as on the level of algebras.
Note that $T^k$ is
connected if and only if the vectors $\{a_1,\ldots,a_n\}$ span
$\mathfrak{t}^d$ over the integers.

Consider the restriction of the action of $T^n$ on $\HH^n$ to the
subgroup $T^k$. This action is hyperhamiltonian with hyperk\"ahler moment map
$$\bar{\mu}_{\R}(z,w) = \iota^* \left(\frac{1}{2} \sum_{i=1}^n (|z_i|^2 - |w_i|^2) u_i
\right) \hspace{15pt}
\text{and}\hspace{15pt}\bar{\mu}_{\C}(z,w) = \iota^* \left(\sum_{i=1}^n (z_i w_i)u_i\right),$$
where $\{u_i\}$ is the standard basis in \(\tnd \cong \R^n.\)
In contrast with the K\"ahler situation, the hyperk\"ahler moment
map is surjective onto $\tnd\oplus(\tn_{\C})^*$.

We denote by 
$M$ the hyperk\"ahler
reduction of $\HH^n$ by the subtorus $T^k$ 
at \((\alpha, 0)
\in (\liet^k)^{*}\oplus (\tk_{\C})^*,\) which is the 
hyperk\"ahler analogue of the K\"ahler toric variety $X = \Cn\mod_{\!\!\a}T^k$. 
Choose a lift $\tilde\a\in\tnd$ of $\a$ along $\iota^*$.
Then $M$ has a natural residual
action of $T^d$ with hyperk\"ahler moment map \(\mhk = \mr\oplus\mc.\)
Note that the choice of subtorus $\Tk\subs\Tn$ is equivalent
to choosing a central arrangement of cooriented hyperplanes in $\tdd$,
where the $i^{\text{th}}$ hyperplane is the annihilator
of $a_i\in\td$.  (The coorientation comes from the fact that we know
for which $x$ we have $\< x,a_i\> > 0$.)
The choice of $\tilde\a$ corresponds to an 
affinization $\arr$ of this arrangement,
where the $i^{\text{th}}$ hyperplane is 
$$H_i = \{x\in\tdd\mid\< x,a_i\>=\< -\tilde\a,\e_i\>\}.$$
Changing $\tilde\a$ by an element $c\in\tdd$ has the effect of translating $\arr$ by $c$,
and adding $c$ to the residual moment map $\mr$.
In order to record the information about coorientations,
we define the half-spaces 
\begin{equation}\label{FG}
F_i = \{x\in\tdd\mid\< x,a_i\>\geq\< -\tilde\a,\e_i\>\}
\hspace{15pt}\text{and}\hspace{15pt} 
G_i = \{x\in\tdd\mid\< x,a_i\>\leq\< -\tilde\a,\e_i\>\},
\end{equation}
which intersect in the hyperplane $H_i$.
Our convention will be to draw pictures, as in Figure 1,
in which we specify the coorientations
of the hyperplanes by shading the polyhedron $\Delta=\cap_{i=1}^n F_i$
(which works as long as $\Delta\neq\emptyset$).
Note that the K\"ahler variety $X$ is precisely the 
K\"ahler toric variety determined by $\Delta$. 

\begin{figure}[h]
\centerline{\epsfig{figure=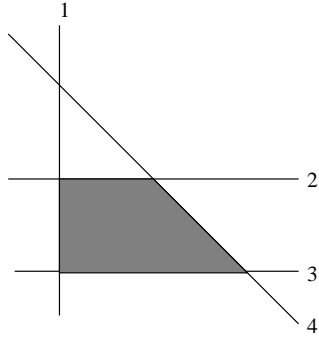}}
\caption{A hypertoric variety of real dimension 8 obtained
by reducing $\H^4$ by $T^2$.}
\end{figure}

The variety $M$ is an orbifold if and only if 
$\arr$ is simple, i.e. if and only if
every subset of $m$ hyperplanes intersect in codimension $m$ \cite[3.2]{BD}.
Furthermore, $M$ is smooth if and only if whenever
some subset of $d$ hyperplanes $\{H_i\}$ has non-empty intersection, 
the corresponding
vectors $\{a_i\}$ form a $\Z$-basis for \(\Z^d \subseteq \liet^d.\)
In this case we will refer to the arrangement itself as {\em smooth}.
We will always assume that $\arr$ is simple, and at times we will also assume that it is smooth.

The hyperplanes $\{H_i\}$
divide \((\liet^d)^{*} \cong \R^d\) into a finite family of closed,
convex polyhedra $$\Delta_A = (\cap_{i\in A}F_i)\cap(\cap_{i\notin A}G_i),$$
indexed by subsets $A\subs\{1,\ldots,n\}$.
Consider the subset $I = 
\{A\subs\{1,\ldots,n\}\mid \Delta_A\text{  bounded}\}$
of the power set of $\{1,\ldots,n\}$.  
For each $A\subs\{1,\ldots,n\}$, let $$M_A = \mr^{-1}(\Delta_A)\cap\mc^{-1}(0).$$
The K\"ahler submanifold
$\left(M_A,\omega_{\R}|_{M_A}\right)$ of $\left(M,\omega_{\R}\right)$
is $d$-dimensional and invariant under the action of $\Td$,
and is therefore $\Td$-equivariantly isomorphic
to the K\"ahler toric variety determined by $\Delta_A$ \cite[6.5]{BD}.
We define the {\em core} $C$ and {\em extended core} $D$
of a hypertoric variety by setting $$C = \cup_{A\in I}M_A\hs\hs\text{ and }\hs\hs
D = \cup_{A}{M_A} = \mc^{-1}(0) =
\{[z,w]\mid z_iw_i=0\text{  for all  }i\},$$
where $[z,w]$
denotes the $T^k$-equivalence class in $M$ of the element
$(z,w)\in \bar{\mu}_{\text{HK}}^{-1}(\a,0)$.
Bielawski and Dancer \cite{BD} show that
$C$ and $D$ are each $\Td$-equivariant
deformation retracts of $M$.  See Corollary \ref{unstable} for a Morse theoretic proof.

We take a minute to discuss the differences between
the combinatorial data
determining a toric variety $X = \Cn\mod_{\!\!\a}\Tk$ and its hypertoric
analogue $M = \Hn\mmod_{\!\!(\a,0)}\Tk$.
Each is determined by $\arr$,
a simple, cooriented, affine arrangement of $n$ hyperplanes in $\tdd$,
defined up to simultaneous translation.
The toric variety $X$ is in fact determined by less information than
this; it depends only on the polyhedron $\Delta = \cap_{i=1}^n F_i$.
Thus if the last hyperplane $H_n$ has the property that
$\cap_{i=1}^{n-1}F_i\subs F_n$, then this hyperplane is superfluous to $X$.
This is not the case for $M$, which means that it is slightly
misleading to call $M$ the hyperk\"ahler analogue of $X$;
more precisely, it is the hyperk\"ahler analogue of
{\em a given presentation} of $X$ as a K\"ahler reduction of $\Cn$.
On the other hand, the $\Td$-equivariant diffeomorphism type of $M$ also 
does not depend on all of the information of $\arr$, as evidenced by the two following results.

\begin{lemma}\label{placement}
The hypertoric varieties \(M_{\a} = \Hn \mmod_{\!\!(\a,0)}\Tk\) 
and $M_{\a'}= \Hn \mmod_{\!\!(\a',0)}\Tk$ are $\Td$-equivariantly
diffeomorphic, and their
cohomology rings can be naturally identified.
\end{lemma}

\begin{lemma}\label{coorientation}
The hypertoric variety $M$ does not depend on the coorientations
of the hyperplanes $\{H_i\}$.
\end{lemma}

This means that, unlike that of $X$, 
the $\Td$-equivariant diffeomorphism type of $M$ depends only on the unoriented
central arrangement underlying $\arr$.
A weaker version of Lemma \ref{placement}, involving
the (nonequivariant) homeomorphism type of $M$, appears in \cite{BD}.

\begin{proofp}
The set of nonregular values for \(\bar{\mu}_{\text{HK}}\) 
has codimension 3 inside of \(\tdd
\oplus (\td_{\C})^*.\) This tells us that the set of regular values 
is simply connected, and we can
choose a path connecting any two regular values \((\alpha,0)\)
and \((\alpha',0),\) 
unique up to homotopy.

Since the moment map \(\bar{\mu}_{\text{HK}}\)
 is not proper, we must take some care in showing that
two fibers are diffeomorphic. To this end, we note that the
norm-square function \(\psi(z,w) = \|z\|^2 + \|w\|^2\) is $T^n$-invariant
and proper on $\Hn$. 
Let $\H_{reg}^n$ denote the open submanifold of $\Hn$ consisting of the
preimages of the regular values of \(\bar{\mu}_{\text{HK}}.\) 
By a direct computation, it is easy
to see that the kernels of $d\psi$ and \(d(\bar{\mu}_{\text{HK}})\) intersect
transversely at any point \(p \in \H_{reg}^n.\) Using the standard
$T^n$-invariant metric on $\Hn$, we define an Ehresmann connection
on $\H_{reg}^n$ with respect to \(\bar{\mu}_{\text{HK}}\) 
such that the horizontal subspaces
are contained in the kernel of $d\psi$.

This connection allows us to lift a path connecting the two
regular values to a horizontal vector field on its preimage in
$\H_{reg}^n$. Since the horizontal subspaces are tangent to
the kernel of $d\psi$, the flow preserves level sets of $\psi$. Note that the
function \[\bar{\mu}_{\text{HK}} \oplus \psi: \Hn \to \tdd \oplus (\td_{\C})^{*} \oplus
\R\] {\it is} proper. By a theorem of Ehresmann \cite[8.12]{BJ}, the
properness of this map
implies that the flow of this vector field exists for all time, and
identifies the inverse image of $(\alpha,0)$ with that of
$(\alpha',0)$. Since the metric, $\psi$, and \(\bar{\mu}_{\text{HK}}\) are
all $T^n$-invariant, the Ehresmann connection is also $T^n$-invariant,
therefore the diffeomorphism identifying the fibers
is $T^n$-equivariant, making 
the reduced spaces are $T^d$-equivariantly diffeomorphic. 
\end{proofp}

\vspace{-\baselineskip}
\begin{proofc}
It suffices to consider the case when we change the orientation of a
single hyperplane within the arrangement.  Changing the coorientation
of a hyperplane $H_l$ is equivalent to defining a new map
$\beta':\tn\to\td$, with $\beta'(\e_i)=a_i$ for $i\neq l$, and $-a_i$ for $i=l$.
This map exponentiates to a map $\beta':\Tn\to\Td$,
and we want to show that the hyperk\"ahler variety obtained
by reducing $\Hn$ by the torus $\ker(\beta')$ is isomorphic to $M$,
which is obtained by reducing $\Hn$ by the torus $\Tk=\ker(\beta)$.
To see this, note that $\ker(\beta')$ and $\ker(\beta)$
are conjugate inside of $\operatorname M(n,\H)$ by the element
$(1,\ldots,1,j,1,\ldots,1)\in \operatorname 
M(1,\H)^n\subs \operatorname M(n,\H)$,
where the $j$ appears in the $l^{\text{th}}$ slot. 
\end{proofc}

\vspace{-\baselineskip}
\begin{example}\label{same}
The three cooriented arrangements of Figure 2 all specify
the same hyperk\"ahler variety $M$ up to equivariant
diffeomorphism.  The first has $X\cong \mathbb{F}_1$ (the
first Hirzebruch surface) and the second and the third have $X\cong\C P^2$.
Note that if we flipped the coorientation of $H_3$ in Figure 2(a)
or 2(c), then
we would get a non-compact
$X\cong\tilde{\C}^2$, the blow-up of $\C^2$
at a point.  If we flipped the coorientation of $H_3$ in
Figure 2(b), then $X$ would be empty.
We make no assumptions about $X$ in this section.

\begin{figure}[h]
\centerline{\epsfig{figure=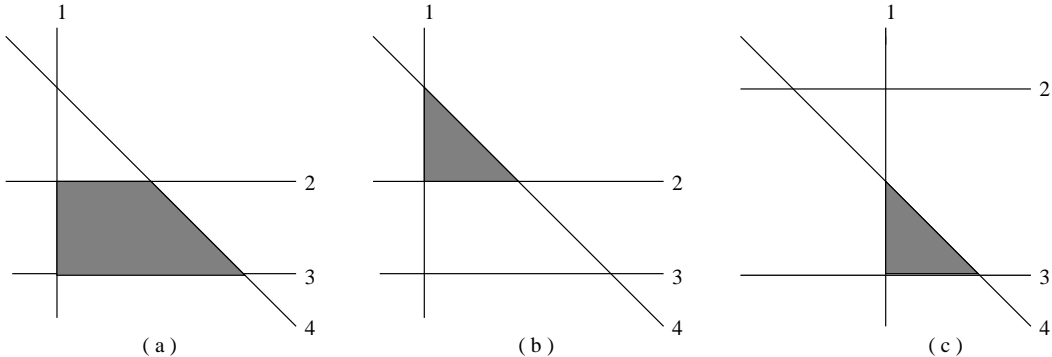}}
\caption{Three arrangements
related by flipping coorientations and translating hyperplanes.}
\end{figure}
\end{example}

The purpose of this paper is to study not just the topology
of $M$, but the topology of $M$ along with the natural hamiltonian
$S^1$ action defined in Section \ref{reduction}.
In order to define this $S^1$ action, it is necessary that
we reduce at a regular value of the form $(\a,0)\in\tdd\oplus(\liet_{\C}^d)^*$,
and although the set of regular values of $\bar{\mu}_{\text{HK}}$ is simply connected,
the set of regular values of the form $(\a,0)$ is not even connected.
Furthermore, 
left multiplication by the diagonal matrix $(1,\ldots,1,j,1,\ldots,1) \in 
\operatorname U(n,\H)$ is not
an $S^1$-equivariant automorphism of $\Hn$, therefore the geometric
structure of $M$ along with a circle action may depend nontrivially
both on the affine structure
and the coorientations of the arrangement $\arr$.
Indeed it must, because we can recover $X$ from $M$ by taking the minimum
$\Phi^{-1}(0)$ of the $S^1$ moment map $\Phi:M\to\R$.
In this sense, the structure of a hypertoric variety $M$
along with a circle action is the universal geometric object from which
both $M$ and $X$ can be recovered.
\end{section}

\begin{section}{Gradient flow on the core}\label{coreflow}
Although $S^1$ does not act on $M$ as a subtorus of $\Td$, we show
below that when restricted to any single 
component $M_A$ of the extended core, $S^1$
{\em does} act as a subtorus of $\Td$, with the subtorus depending
combinatorially on $A$.  This will allow us to give a combinatorial
analysis of the gradient flow of $\Phi$ on the extended core.

\begin{lemma}\label{sides}
Let $x$ be an element of $\tdd$, and consider a point $[z,w]\in\mr^{-1}(x)\cap D$.
We have $x \in F_i$ if and only if $w_i=0$, and $x\in G_i$
if and only if $z_i=0$.
\end{lemma}

\begin{proof}
The fact that $[z,w]\in D$ tells us that $z_iw_i=0$.
Then 
\begin{eqnarray*}
x\in F_i & \iff & \< -\tilde\a,\e_i\>\leq\<\mr[z,w],
	           a_i\>  =\<\bar\mr(z,w)-\tilde\a,\e_i\>  \\ 
        & \iff  &  \half|z_i|^2-\half|w_i|^2=\<\bar\mr(z,w),
		    \e_i\>\geq 0.                                  \\
\end{eqnarray*}
Since $z_iw_i=0$, this is equivalent to the condition $w_i=0$.
The second half of the lemma follows similarly.
\end{proof}

On the suborbifold $M_A\subs D\subs M$ we have $z_i = 0$ for all $i\in A$
and $w_i = 0$ for all $i\notin A$, therefore for $\tau\in S^1$
and $[z,w]\in M_A$,
$$\tau[z,w]=[z,\tau w] = 
[\tau_1 z_1,\ldots \tau_n z_n, \tau_1^{-1}w_1,\ldots \tau_n^{-1}w_n], \text{ where }
\tau_i =
\begin{cases}
\tau^{-1} & \text{if }i\in A,\\
1 & \text{if }i\notin A.
\end{cases}
$$
In other words, the $S^1$ action on $M_A$
is given by the one dimensional subtorus 
$(\tau_1, \ldots, \tau_n)$
of the original torus $T^n$,
hence the moment map $\Phi|_{M_A}$ is given (up to
an additive constant) by 
$$\Phi[z,w] = \left<\mr[z,w],\sum_{i\in A}a_i\right>.$$
This formula allows us to compute the fixed point sets of
the $S^1$ action.
Since $S^1$ acts freely on \((\liet_{\C}^d)^{*} \backslash \{0\}\)
and \(\mc: M \to (\liet_{\C}^d)^{*}\) is $S^1$-equivariant,
we must have $M^{S^1} \subseteq \mc^{-1}(0) = D$.
For any subset $B\subs\{1,\ldots,n\}$, let $M_A^B$ be the toric subvariety
of $M_A$ defined by the conditions $z_i=w_i=0$ for all $i\in B$.
Geometrically, $M_A^B$ is defined by the (possibly empty) face
$\cap_{i\in B}H_i \cap \Delta_A$ of the polyhedron $\Delta_A$.

\begin{proposition}\label{fixed}
The fixed point set of the action of $S^1$ on $M_A$
is the union of those toric subvarieties $M_A^B$ such that
$\sum_{i\in A}a_i\in \liet^d_B:=\operatorname{Span}_{j\in B}a_j$.
\end{proposition}

\begin{proof}
The moment map $\Phi|_{M_A^B}$ will be constant if and only if
$\sum_{i\in A}a_i$ is perpendicular to 
$\ker\left(\tdd\onto(\liet^d_B)^*\right)$,
i.e. if $\sum_{i\in A}a_i$ lies in the kernel of the projection
$\td\onto\td/\liet^d_B$.
\end{proof}

\vspace{-\baselineskip}
\begin{corollary}
Every vertex $v\in\tdd$ of the polyhedral complex defined by $\arr$
is the image of an $S^1$-fixed point in $M$.
Every component of $M^{S^1}$ has dimension less than or equal to $d$,
and the only component of dimension $d$ is $M_{\emptyset}=X=\Phi^{-1}(0)$. 
\end{corollary}

For any point $p\in M^{S^1}$, the {\em stable orbifold} 
$S(p)$ at $p$ is defined to be the set of $x\in M$ such
that $x$ approaches $p$ when flowing along the vector field 
$-\operatorname{grad}(\Phi)$,
and the {\em unstable orbifold} $U(p)$ at $p$ is defined to be the stable orbifold with respect
to the function $-\Phi$.
For any suborbifold $Y\subs M^{S^1}$, the unstable orbifold $U(Y)$ at $Y$ is defined
to be the union of $U(y)$ for all $y\in Y$.  In general, for $y\in Y$, we have the identity
$\dim_{\R} U(Y) + \dim_{\R} S(y) = 4d$.

Let $Y\subs M^{S^1}$ be a component of the fixed point set of $M$.
Let $v\in\tdd$ be a vertex in the polehedron $\mr(Y)$, and let $y$
be the unique preimage of $v$ in $Y$.

\begin{proposition}
The unstable orbifold $U(Y)$ is a complex suborbifold of complex 
dimension at most $d$, contained in the core $C\subs M$.
If $\arr$ is smooth at $y$, then $\dim_{\C}U(Y)=d$, and the closure of $U(Y)$
is an irreducible component of $C$.
\end{proposition}

\begin{proof}
For simplicity, we will assume that $v=\cap_{j=1}^dH_j$.
For all $l\in\otd$, let $b_l\in\td_{\Z}$ be the smallest integer vector
such that $\left<a_j,b_l\>=0$ for $j\neq l$ and $\left<a_l,b_l\> >0$.
Geometrically, $b_l$ is the primitive integer vector on the line
$\cap_{j\neq l}H_j$ pointing in the direction of $\Delta$.
Note that $M$ is smooth at the $T^d$-fixed point above $v$
if and only if $\left<a_l,b_l\>=1$ for all $l\in\otd$.
Let $R_l\subs\tdd$ be the ray eminating from $v$ in the direction of $b_l$,
and ending before it hits another vertex.  Let $Q_l$ be the analogous ray in the opposite
direction.

Let $\Delta_A$ be a region (not necessarily bounded) of the polyhedral 
complex defined by $\arr$ adjacent to $R_l$.
The preimage $\mr^{-1}(R_l)\cap D$ of $R_l$ in $D$ is a complex line,
and it is contained in the unstable
orbifold at $U(Y)$ if and only if $\left<b_l,\sum_{i\in A}a_i\right>\geq 0$.
If $\left<b_l,\sum_{i\in A}a_i\right> <0$, it is contained in the stable orbifold $S(y)$.
The preimage $\mr^{-1}(Q_l)\cap D$ of $Q_l$ in $D$ is also a complex line,
contained in the unstable
orbifold $U(Y)$ if and only if $\<-b_l,a_l+\sum_{i\in A}a_i\>\geq 0$,
and otherwise in $S(y)$.
Since $\<b_l,\sum_{i\in A}a_i\>+\<-b_l,a_l+\sum_{i\in A}a_i\>
=-\left<a_l,b_l\> <0$, at most one of these two directions can be unstable.
In the smooth case, $\left<a_l,b_l\>=1$ for all $l$, and exactly one of
the two directions is unstable.

Consider the polytope $\Delta_v$ incident to $v$ and characterized by the property that its
edges at the vertex $v$ are exactly the unstable directions. 
The toric variety $X_{\Delta_v}\subs D$ is contained
in the closure of $U(Y)$, and a dimension count tells us that this containment is an equality.
In the smooth case, $\Delta_v$ is $d$-dimensional, and $X_{\Delta_v}$ is a component of the core.
\end{proof}

Note that, even in the smooth case, it is not necessarily the case that the
$R_l$ direction is stable and the $Q_l$ direction is unstable.  See, for example,
the vertex $v = H_1\cap H_2$ in Figure 2(c).

\begin{corollary}
There is a natural injection from the set of bounded regions 
$\{\Delta_A\mid A\in I\}$
to the set of connected components of $M^{S^1}$.
If $\arr$ is smooth, this map is a bijection.
\end{corollary}

\begin{proof}
To each $A\in I$, we associate the fixed subvariety $M_A^B$ corresponding
to the face of $\Delta_A$ on which the linear functional $\sum_{i\in A}a_i$
is minimized, so that $M_A$ is the closure of
$U(M_A^B)$.  If $\arr$ is smooth, then every connected component of the fixed point
set will have a component of the core as its closed unstable orbifold.
\end{proof}

\vspace{-\baselineskip}
\begin{corollary}\label{unstable}
The core of $M$ is equal to the union of the unstable orbifolds of
the connected components of $M^{S^1}$, hence $C$ is a $\Tdso$-equivariant
deformation retract of $M$.
\end{corollary}

\begin{example}
In Figure 3, representing a reduction of $\H^5$ by $T^3$,
we choose a metric on $(\liet^2)^*$ in order
to draw the linear functional $\sum_{i\in A}a_i$ as a vector
in each region $\Delta_A$.  We see that $M^{S^1}$ has three components,
one of them $X\cong\mathbb F_1$, one of them a projective line, with another $\mathbb F_1$
as its unstable manifold, and one of them a point, with a $\C P^2$ as its unstable manifold.

\begin{figure}[h]
\centerline{\epsfig{figure=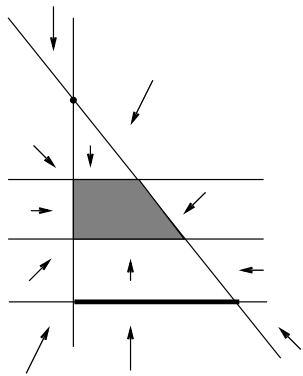}}
\caption{The gradient flow of $\Phi:M\to\R$.}
\end{figure}
\end{example}

\begin{example}
The hypertoric variety represented by Figure 4 has a fixed point set with four
connected components (three points and a $\C P^2$), but only three components in its core.
This phenomenon can be blamed on the orbifold point represented by the
intersection of $H_3$ and $H_4$, which has a one-dimensional unstable orbifold.

\begin{figure}[h]
\centerline{\epsfig{figure=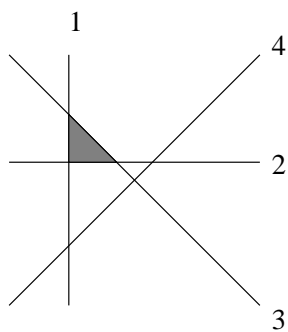}}
\caption{A singular example.}
\end{figure}
\end{example}
\end{section}

\begin{section}{Equivariant cohomology}\label{Konno}
In this section we extend Konno's computations
of the ordinary and $\Td$-equivariant cohomologies of $M$
to the $S^1$-equivariant setting.
We follow Konno's approach of restricting to the smooth case
to simplify arguments involving line bundles on $M$.
Hausel and Sturmfels, however, prove
theorems analogous to \ref{htm} and \ref{hm} with rational coefficients 
in the orbifold case, and Theorems \ref{htsm} and \ref{hsm}
extend to this setting as well (see Remark \ref{lawrence}).

\begin{theorem}\label{htm}{\em\cite{K2}}
The $\Td$-equivariant cohomology ring of a smooth 
hypertoric variety $M$ is given by
$$\htm = \Z[u_1,\ldots,u_n]\bigg/
\left<\prod_{i\in S}u_i\hs\bigg{\vert}\hs\bigcap_{i\in S} H_i = \emptyset \right>.$$
\end{theorem}

\begin{remark}
This is precisely
the Stanley-Reisner ring of
the unoriented matroid determined by the arrangement $\arr$ \cite{HS}.
\end{remark}

Just as the cohomology of a toric variety is obtained
from the equivariant cohomology by introducing linear relations
that generate $\ker\iota^*=\left(\ker\beta\right)^{\bot}$, the same
is true for hypertoric varieties:

\begin{theorem}\label{hm}{\em\cite{K1}}
The ordinary cohomology ring of a smooth hypertoric variety $M$ is given by
$$\hm = \htm \big/ \<\Sigma a_iu_i\in \ker\iota^*\>.$$
\end{theorem}

The rest of this section will be devoted to the proof of the following two theorems.

\begin{theorem}\label{htsm}
Let $M$ be the hypertoric variety corresponding to a
smooth, cooriented arrangement $\arr$.
Given any minimal set $S\subs \{1,\ldots,n\}$ such that 
$\cap_{i\in S} H_i = \emptyset$, let $S = S_1\sqcup S_2$
be the unique splitting of $S$ such that
$\big(\cap_{i\in S_1}G_i\big)\cap
\big(\cap_{j\in S_2}F_j\big) = \emptyset$ (see \eqref{FG}).
Then the $\Tdso$-equivariant cohomology of $M$ is given by
$$\htsm\cong\Z[u_1,\ldots,u_n,x]\bigg/
\< \prod_{i\in S_1}u_i\times\prod_{j\in S_2}(x-u_j)
\hs\bigg{\vert}\hs\bigcap_{i\in S} H_i = \emptyset\>.$$
\end{theorem}

\begin{theorem}\label{hsm}
In the notation of Theorem \ref{htsm},
the $S^1$-equivariant cohomology ring
of $M$ is given by
$$\hsm \cong \htsm \big/ \<\Sigma a_iu_i\in \ker\iota^*\>.$$
\end{theorem}

\begin{remark}
Konno observes that the quotient map from the abstract
polynomial ring $\Z[u_1,\ldots,u_n]\to\htm$ is precisely the $\Td$-equivariant
Kirwan map $$\ktd:H^*_{\Tn}(\cot)\to\htm$$
induced by the inclusion
$\mu^{-1}(\a,0)\hookto\cot$.  
Likewise, the map from $\Z[u_1,\ldots,u_n]/\ker\iota^*$ to $H^*(M)$
is the ordinary Kirwan map
$$\k:H^*_{\Tk}(\cot)\to H^*(M).$$
The analogous maps for K\"ahler reductions
are known to always be surjective \cite[5.4]{Ki}, but the hyperk\"ahler case
remains open.  Thus Theorems \ref{htm} and \ref{hm} can be interpreted as saying
that the Kirwan maps for hypertoric varieties are surjective, and computing the kernel.
Likewise, Theorems \ref{htsm} and \ref{hsm} assert that the $S^1$-equivariant Kirwan maps
$$\ktdso:H^*_{\Tn\times S^1}(\cot)\to\htsm$$ and
$$\kso:H^*_{\Tk\times S^1}(\cot)\to\hsm$$ are surjective,
and provide computations of their kernels.
\end{remark}

In order to apply Konno's results, we will make use of
the principle of {\em equivariant formality},
proven for compact manifolds in \cite{Ki}, which we adapt to our situation
in Proposition \ref{formality}.
For the sake of simplicity, we will restrict our attention
to the case where $X$ is compact and nonempty.
This condition will be necessary for the application of
Proposition \ref{formality} and the proof of Theorem \ref{tara},
both of which require a proper Morse function, which we get from
Proposition \ref{proper}.
We note, however, that both Proposition \ref{formality} and Theorem
\ref{tara} can be extended to the case of a general hypertoric
variety by a Mayer-Vietoris argument, using the fact that the core
$C\subs M$ is a compact $\Td\times S^1$-equivariant deformation retract.
We present the slightly less general Morse theoretic proofs
only because we find them more pleasant.

\begin{proposition}\label{formality}
Let $M$ be a symplectic
orbifold, possibly noncompact but of finite topological type.
Suppose that $M$ admits a hamiltonian action of a torus $T\times S^1$,
and that the $S^1$-component $\Phi:M\to\R$ of the moment map
is proper and bounded below.  Then 
$H^*_{T\times S^1}(M)$ 
is a free module over $H^*_{S^1}(pt)$.
\end{proposition} 

\begin{proof}
Because $\Phi$ is a moment map, it is a Morse-Bott function
such that all of the critical suborbifold and their normal bundles
carry almost complex structures.  Thus we get a Morse-Bott stratification
of $M$ into even-dimensional $T$-invariant suborbifolds.
This tells us, as in [Ki, 5.8],
that the spectral sequence associated to the fibration
$M\hookto EG\times_G M\to BG$ collapses, and we get the desired result.
\end{proof}

Consider the following commuting square of maps, where
$\phi$ and $\psi$ are each given by setting $x$ to zero.

$$\begin{CD}
H^*_{\Tn\times S^1}(\cot) @>\ktdso >> \htsm\\
@V\phi VV @ VV\psi V\\
H^*_{\Tn}(\cot) @>\ktd >> \htm\\
\end{CD}$$

Proposition \ref{formality} has the following consequence.

\begin{corollary}\label{enough}
Let $\mathcal{I}\subs\ker\ktdso$ be an ideal with $\phi(\mathcal{I}) = \ker\ktd$.
Then $\mathcal{I}=\ker\ktdso$.
\end{corollary}

\begin{proof}
Suppose that $a\in\ker\ktdso\smallsetminus\mathcal{I}$ is a homogeneous class of minimal degree,
and choose $b\in\mathcal{I}$ such that $\phi(a-b)=0$.  Then $a-b = cx$ for some
$c\in H^*_{\Tn\times S^1}(\cot)$.  By Proposition \ref{formality},
$cx\in\ker\ktdso\impl c\in\ker\ktdso$, hence $c\in\ker\ktdso\smallsetminus\mathcal{I}$
is a class of lower degree than $a$.
\end{proof}

\vspace{-\baselineskip}
\begin{lemma}\label{surjective}
The equivariant Kirwan map $\ktdso$ is surjective.
\end{lemma}

\begin{proof}
Suppose that $\gamma\in\htsm$ is a homogeneous class of minimal degree
that is {\em not}
in the image of $\ktdso$.  By Theorem \ref{htm} $\ktd$ is surjective, hence we may choose a class
$\eta\in\phi^{-1}\ktd^{-1}\psi(\gamma)$.
Then $\ktdso(\eta)-\gamma = x\delta$ for some $\delta\in\htsm$,
and therefore $\delta$ is a class of lower degree that is not in the image of $\ktdso$.
\end{proof}

\vspace{-\baselineskip}
\begin{proofhtsm}
For any element $h\in H^2_{T^n\times S^1}(\Hn;\Z)$,
let $\tilde{L}_h = \Hn\times\C_{h}$ 
be the $\Tn\times S^1$-equivariant line bundle on $\Hn$
with equivariant Euler class $h$. 
This gives 
$\tilde{L}_h$, as well as its dual $\tilde{L}_h^*$, 
the structure of a \(T^n \times S^1\)-equivariant bundle. 
Let $$L_h = \tilde{L}_h|_{\bar{\mu}_{\C}^{-1}(0)^{\text{ss}}}/T_{\C}^k$$ 
be the quotient $\Td\times S^1$-equivariant line bundle on $M$.
Let $\{u_i\}$ be the standard basis of $\Tnd$.
Identifying $H^2_{T^n\times S^1}(\Hn;\Z)$ with $\Tnd\oplus\Z x$,
We will use $\tilde L_i$ to denote the bundle $\tilde L_{u_i\oplus 0}$, 
and $\tilde K$ to denote the bundle $\tilde L_{0\oplus x}$,
with quotients $L_i$ and $K$. 
Since the $\Tdso$-equivariant Euler class $e(L_i)$
is the image of $u_i\oplus 0$ under the hyperk\"ahler Kirwan map
$H^*_{T^n\times S^1}(\Hn)\to \htsm$, we will abuse notation and 
denote it by $u_i$.  Similarly, we will denote $e(K)$ by $x$.
Corollary \ref{surjective} tells us that
$\htsm$ is generated by $u_1,\ldots,u_n,x$.

Consider the $\Tn\times S^1$-equivariant section $\tilde s_i$ of $\tilde{L}_i$
given by the function $\tilde s_i(z,w)=z_i$.
This descends to a $\Tdso$-equivariant section $s_i$ of $L_i$
with zero-set $Z_i
:= \{[z,w]\in M\mid z_i=0\}$. 
Similarly, the function $\tilde t_i(z,w)=w_i$ defines a $\Tdso$-equivariant
section of $L_i^*\otimes K$ with zero set $W_i := \{[z,w]\in M\mid w_i=0\}$.
Thus the divisor $Z_i$ represents the cohomology class $u_i$,
and $W_i$ represents $x-u_i$.
Note, by the proof of Lemma \ref{sides},
that $\mr(Z_i) = G_i$ and $\mr(W_i) = F_i$ for all $1\leq i \leq n$.

Let $S = S_1\sqcup S_2$ be a subset of $\{1,\ldots n\}$ such that
$\big(\cap_{i\in S_1}G_i\big)\cap
\big(\cap_{j\in S_2}F_j\big) = \emptyset,$ and hence
$$\big(\cap_{i\in S_1}Z_i\big)\cap
\big(\cap_{j\in S_2}W_j\big)\subs\mr^{-1}\bigg(
\big(\cap_{i\in S_1}G_i\big)\cap
\big(\cap_{j\in S_2}F_j\big)\bigg) = \emptyset.$$
Now consider the vector bundle
$E_S=\left(\oplus_{i\in S_1}L_i\right)\oplus\left(\oplus_{j\in S_2}L_j^*\otimes K\right)$
with equivariant Euler class $$e(E_S)=\prod_{i\in S_1}u_i\times\prod_{j\in S_2}(x-u_j).$$ 
The section $\left(\oplus_{i\in S_1}s_i\right)\oplus\left(\oplus_{i\in S_2}t_i\right)$ 
is a nonvanishing equivariant
global section of $E_S$,
hence for any such $S$, $e(E_S)$ is trivial in $\htsm$.

The fact that $u_1,\ldots,u_n,x$ generate $\htsm$ is proven in Lemma \ref{surjective},
and the fact that
we have found all of the relations
follows from Theorem \ref{htm} and Corollary \ref{enough}.
\end{proofhtsm}

\vspace{-\baselineskip}
\begin{proofhsm}
The proof of this theorem is identical to the proof of Theorem \ref{htsm},
making use of Theorem \ref{hm} rather than Theorem \ref{htm}.
\end{proofhsm}

How sensitive are the invariants $\htsm$ and $\hsm$?
We can recover $\htm$ and $\hm$
by setting $x$ to zero, hence they are at least as fine
as the ordinary or $\Td$-equivariant cohomology rings.
The ring $\htsm$ does {\em not} depend on coorientations,
for if $M'$ is related to $M$ by flipping the coorientation
of the $l^{\text{th}}$ hyperplane $H_k$, then the map taking $u_i$ to
$u_i$ for $i\neq l$ and $u_l$ to $x-u_l$ is an isomorphism
between $\htsm$ and $H^*_{\Td\times S^1}(M')$.
It is, however, dependent on the affine structure of 
the arrangement $\arr$.

\begin{example}\label{ts}
We compute the equivariant cohomology ring $\htsm$
for the hypertoric varieties $M_a$, $M_b$, and $M_c$
defined by the arrangements in Figure 2(a), (b), and (c), respectively.
$$H^*_{\Td\times S^1}(M_a)=\Z[u_1,\ldots,u_4,x]\big/
\< u_2u_3, u_1(x-u_2)u_4, u_1u_3u_4 \>,$$
$$H^*_{\Td\times S^1}(M_b)=\Z[u_1,\ldots,u_4,x]\big/
\< (x-u_2)u_3, u_1u_2u_4, u_1u_3u_4 \>,$$
$$H^*_{\Td\times S^1}(M_c)=\Z[u_1,\ldots,u_4,x]\big/
\< u_2u_3, (x-u_1)u_2(x-u_4), u_1u_3u_4 \>.$$
As we have already observed, $H^*_{\Td\times S^1}(M_a)$ and $H^*_{\Td\times S^1}(M_b)$
are isomorphic by interchanging $u_2$ with $x-u_2$.
One can check that the annihilator of $u_2$ in $H^*_{\Td\times S^1}(M_a)$
is the principal ideal generated by $u_3$, while the ring $H^*_{\Td\times S^1}(M_c)$
has no degree $2$ element whose annihilator is generated by a single element of degree $2$.
Hence $H^*_{\Td\times S^1}(M_c)$ is not isomorphic to the other two rings.
\end{example}

The ring $\hsm$, on the other hand, is sensitive to coorientations
as well as the affine structure of $\arr$.

\begin{example}\label{s}
We now compute the ring $\hsm$ for 
$M_a$, $M_b$, and $M_c$ of Figure 2.
Theorem \ref{hsm} tells us that we need only to quotient the ring
$\htsm$ by $\ker(\iota^*)$.
For $M_a$, the kernel of $\iota_a^*$ is generated by
$u_1+u_2-u_3$ and $u_1-u_4$, hence we have
\begin{eqnarray*}
H^*_{S^1}(M_a)&=&\Z[u_2,u_3,x]\big/
\< u_2u_3, (u_3-u_2)^2(x-u_2), (u_3-u_2)^2u_3 \>\\
&\cong & \Z[u_2,u_3,x]\big/\< u_2u_3, (u_3-u_2)^2(x-u_2), u_3^3 \>.\\
\end{eqnarray*}
Since the hyperplanes of 2(c) have the same coorientations as those
of 2(a), we have $\ker\iota_b^*=\ker\iota_a^*$, hence
\begin{eqnarray*}
H^*_{S^1}(M_c) &=& \Z[u_2,u_3,x]\big/
\< u_2u_3, (x-u_3+u_2)^2u_2, (u_3-u_2)^2u_3 \>\\
&\cong & \Z[u_2,u_3,x]\big/\< u_2u_3, (x-u_3+u_2)^2u_2, u_3^3 \>.\\
\end{eqnarray*}
Finally, since Figure 2(b) is obtained from 2(a) by flipping the coorientation
of $H_2$, we find that $\ker(\iota_b^*)$ is generated by $u_1-u_2-u_3$
and $u_1-u_4$, therefore
$$H^*_{S^1}(M_b)=\Z[u_2,u_3,x]\big/
\< (x-u_2)u_3, (u_2+u_3)^2u_2, (u_2+u_3)^2u_3 \>.$$
As in Example \ref{ts}, $H^*_{S^1}(M_a)$ and $H^*_{S^1}(M_c)$
can be distinguished by the fact that the annihilator
of $u_2\in H^*_{S^1}(M_a)$ is generated by a single element of degree $2$,
and no element of $H^*_{S^1}(M_c)$ has this property.
On the other hand, $H^*_{S^1}(M_b)$ is distinguished from 
$H^*_{S^1}(M_a)$ and $H^*_{S^1}(M_c)$ by the fact that neither $x-u_2$ nor $u_3$ cubes to zero.
\end{example}

\begin{remark}\label{lawrence}
Theorems \ref{htsm} and \ref{hsm} can be interpreted
in light of the recent work of Hausel and Sturmfels \cite{HS} 
on Lawrence toric varieties.
The Lawrence toric variety $N$ associated to the arrangement $\arr$
is the K\"ahler reduction $\cot\mod\Tk$,
so that $M$ sits inside of $N$ as the complete
intersection cut out by the equation $\bar\mc(z,w)=0$.
The residual torus acting on $N$ has dimension $d+n$,
and includes the $(d+1)$-dimensional torus $\Td\times S^1$
acting on $M$, and
the inclusion of $M$ into
$N$ induces an isomorphism on $\Td\times S^1$-equivariant cohomology.
One can use geometric arguments similar to those that were applied to prove
Theorem \ref{htsm},
or the purely combinatorial approach of \cite{HS},
to show that
$$H_{T^{d+n}}^*(N)=\Q[u_1,\ldots,u_n,v_1,\ldots,v_n]\bigg/
\< \prod_{i\in S_1}u_i\times\prod_{j\in S_2}v_j\hs\bigg{\vert}\hs\bigcap_{i\in S} H_i = \emptyset\>.$$
From here we can recover $\htsm = H_{\Td\times S^1}^*(N)$
by setting $u_i+v_i = u_j+v_j$ for all $i,j\leq n$.
Note that Hausel and Sturmfels' work applies to the general orbifold case.
\end{remark}
\end{section}

\begin{section}{A deformation of the Orlik-Solomon algebra of \boldmath$\arr$}\label{defos}
Let $M_{\R} \subs M$ be the {\em real locus}
$\{[z,w]\in M\mid z,w \text{ real}\}$ of $M$ with respect to the complex structure $J_1$.
The full group $\Td\times S^1$ does not act on $M_{\R} $,
but the subgroup $\Tdr\times\Zt$ does act,
where $\Tdr :=\Z_2^d\subs \Td$ is the fixed point set
of the involution of $\Td$ given by complex
conjugation.\footnote{It is interesting to note that
the real locus with respect to the complex structure $J_1$
is in fact a complex submanifold with respect to the
one of the other complex structures on $M$.  
The action of $\Tdr$ is holomorphic because $\Tdr$ is a subgroup of $\Td$,
which preserves all of the complex structures on $M$.
The action of $\Zt$, on the other hand, is anti-holomorphic,
i.e. it can be thought of as complex conjugation.}
In this section we will study the geometry of the real locus,
focusing in particular on the properties of the residual $\Zt$ action.

A proof of a more general statement of the following theorem is 
forthcoming in \cite{HH}. 


\begin{theorem}\label{tara}
Let $G = \Tdso$ or $\Td$,
and $G_{\R} = \Tdr\times\Zt$ or $\Tdr$.
Then we have $H^{2*}_{G}(M;\Zt)\cong
H^{*}_{G_{\R}}(M_{\R} ;\Zt),$
i.e. the rings are isomorphic
by an isomorphism that halves the grading.
Furthermore, this isomorphism identifies the class 
$u_i\in H^{*}_{G}(M;\Zt)$, 
represented by the divisor $Z_i$, with the class in $H^*_{G_{\R}}(M_{\R} ;\Zt)$
represented by the divisor $Z_i\cap M_{\R} $, and likewise takes
$x-u_i$ (if $G=\Tdso$) or
$-u_i$ (if $G=\Td$) to the class represented by $W_i\cap M_{\R} $. 
\end{theorem}

\begin{sketch}
Consider the injection $H^*_G(M;\Zt)\hookto H^*_G(M^G;\Zt)$
given by the inclusion of the fixed point set into $M$.
The essential idea is to show that 
a class in $H^*_G(M^G;\Zt)$ extends over $M$ if and only
if it extends to the set of points on which $G$ acts with
a stabilizer of codimension at most $1$,
and then to show that a similar
statement in $G_{\R}$-equivariant cohomology 
also holds for the real locus $M_{\R} $ with its $G_{\R}$
action. One then uses a canonical isomorphism \(H^2_G(pt,\Z_2) \cong 
H^1_{G_{\R}}(pt,\Z_2)\) to give the result. 

The key to the proof is a noncompact $G_{\R}$ version
of the proposition in 
\cite{TW} stating that the $G_{\R}$-equivariant Euler class of 
the negative normal bundle of a critical point $p$ 
is not a zero divisor, which can be shown explicitly
using a local normal form for the actions of $G$ and $G_{\R}$. 
The proposition then follows 
from standard $G_{\R}$ versions of the Thom isomorphism theorem
with coefficients in $\Z_2$. 
Since a component of the 
moment map is proper, bounded below, and has finitely many
fixed points, one can then check that the inductive argument, given 
in Section 3 of \cite{TW}
to complete the proof of \cite[Thm 1]{TW} also holds in this case.
\end{sketch}

%
%
%

Let us consider the restriction of the hyperk\"ahler moment map
$\mhk=\mr\oplus\mc$ to $M_{\R} $.  Since $z$ and $w$ are real
for every $[z,w]\in M_{\R} $, the map $\mc$
takes values in $\td_{\R}\subs\td_{\C}$, which we will identify
with $i\R^n$, so that $f = \mr|_{M_{\R}} \oplus\mc|_{M_{\R}}$ takes values
in $\R^n\oplus i\R^n\cong\C^n$.
Note that $f$ is $\Zt$-equivariant, with $\Zt$ acting on
$\Cn$ by complex conjugation.

\begin{lemma}
The map $f:{M_{\R}}\to\Cn$ is surjective, and the fibers are the orbits of $\Tdr$.
The stabilizer of a point $x\in {M_{\R}}$ has order 
$2^r$, where $r$ is the number of hyperplanes in the complexified
arrangement $\arr_{\C}$ containing the point $f(x)$.
\end{lemma}

\begin{proof}
For any point $p=a+bi\in\Cn$, choose a point $[z,w]\in M$ such that
$\mr[z,w]=a$ and $\mc[z,w]=b$.  We can hit $[z,w]$
with an element of $\Td=\Tn/\Tk$ to make $z$ real, and the fact that
$\mc[z,w]\in\R^n$ forces $w$ to be real as well, hence we may assume
that $[z,w]\in {M_{\R}}$.
Then $f^{-1}(p)=\mu_{\text{HK}}^{-1}\left(a,b\right)\cap {M_{\R}}
= \Td[z,w]\cap {M_{\R}} = \Tdr[z,w]$.
The second statement follows easily from \cite[3.1]{BD}.
\end{proof}

Let $Y\subs {M_{\R}}$ be the locus of points on which $\Tdr$
acts freely, i.e. the preimage under $f$ of the space 
$\comp:=\bomp$.
The inclusion map $Y\hookto {M_{\R}}$
induces maps backward on cohomology, which we will denote
$$\phi:\hr\to H^*_{\Tdr}(Y;\Zt)\cong H^*(\comp;\Zt)$$
and
$$\phi_2:\hrs\to H^*_{\Tdr\times\Zt}(Y;\Zt)\cong H^*_{\Zt}(\comp;\Zt).$$
The ring $H^*(\comp;\Z)$ has been studied extensively, and is called
the {\em Orlik-Solomon algebra} \cite{OT}, which we will denote
by $\os$.
A remarkable fact about the Orlik-Solomon algebra is that it depends
only on the combinatorial structure of $\arr$; the following is a presentation
in terms of {\it anticommuting} generators $e_1,\ldots,e_n$ \cite{OS}:
$$\os\cong H^*(\comp;\Z)\cong\Z[e_1,\ldots,e_n]\big/
\< \Pi_{i\in S}e_i\mid\cap_{i\in S} H_i = \emptyset \>.$$
Rather than working with anticommuting generators, we can work over the ground field $\Zt$,
in which case commutativity and anticommutativity are the same.
Because $\os$ is torsion-free \cite[3.74]{OT}, we have
$$\ost\cong H^*(\comp;\Zt)\cong\Zt[e_1,\ldots,e_n]\big/
\< \Pi_{i\in S}e_i\mid\cap_{i\in S} H_i = \emptyset \>
+\< e_i^2\mid i\leq n\>,$$
where $\operatorname{deg}(e_i)=1$.

\begin{claim}
The map $\phi:H^*_{\Tdr}(M_{\R};\Zt)\to\ost$ 
takes $u_i$ to $e_i$, hence $\ker\phi$
is generated by
the set $\{u_i^2\mid i\leq n\}$.
\end{claim}

\begin{proof}
Recall from Section \ref{ht}
that the hyperplane $H_i\subs \tdd$
is defined by the equation $\< x,a_i\> = \< -\tilde\a,\e_i\>$.
Let $\eta_i:\Cn\to\C$ be the affine map
taking $x$ to $\< x,a_i\> + \< \tilde\a,\e_i\>$,
so that $H_i^{\C}$ is cut out of $\Cn$ by $\eta_i$.
Then $\eta_i$ restricts to a map $\comp\to\C^*$, and
Orlik and Terao identify $e_i$ with the cohomology class
represented by the pull-back
of the submanifold $\R_-$ (the negative reals) along $\eta_i$ [OT, 5.90].
Theorem \ref{tara} tells us that the cohomology
class $u_i$ is represented by the divisor $Z_i\cap {M_{\R}}$.
By Lemma \ref{sides}, $f(Z_i\cap {M_{\R}}) = G_i\cap \R^n$,
hence $\phi(u_i)$ is the class represented by the submanifold
$G_i\cap \R^n\cap \comp = 
\eta_i^{-1}(\R_-)$.
\end{proof}

\vspace{-\baselineskip}
\begin{remark}\label{kernel}
The fact that the classes $u_i^2$ lie in the kernel of $\phi$ can be seen
by noting that $u_i$ is represented
both by the divisor $Z_i = \{[z,w]\in {M_{\R}}\mid z_i=0\}$, 
and, since $u_i = -u_i$ over $\Zt$, 
also by the divisor $W_i = \{[z,w]\in {M_{\R}}\mid w_i=0\}$.
The condition $x\in Z_i\cap W_i$ says exactly that
$\mr(x)\in H_i$, therefore
the intersection $Z_i\cap W_i \cap Y$ is empty.
\end{remark}

In some sense we have cheated here; we have concluded that
we can recover a presentation of $\ost$ from a presentation
of $\htm$, but we used the fact that we already have a presentation
of $\ost$.  In the $\Zt$-equivariant picture, however,
our trivial observation turns magically into new information,
giving us a presentation of the equivariant cohomology
ring $\hscomp$.

\begin{theorem}\label{surj}
The map $\phi_2$ is surjective, with kernel generated by 
$\{u_i (x-u_i)\mid i\leq n\}$.
Hence 
$$\hscomp \cong H^*_{\Tdr\times\Zt}(M_{\R};\Zt)\big/\ker\phi_2$$
$$\cong\Zt[u_1,\ldots,u_n,x]\big/
\<\Pi_{i\in S_1}u_i\times\Pi_{j\in S_2}(x-u_j)
\mid\cap_{i\in S} H_i = \emptyset\>
+\<u_i(x-u_i)\mid i\leq n\>.$$
\end{theorem}

\begin{proof}
The fact that $\arr_{\C}$ is the 
complexification of a real arrangement
tells us that the linear forms $\eta_i$ have real coefficients,
therefore the generators $e_i$ of $\hcomp$ are $\Zt$-invariant.
This implies, by \cite[3.5]{Bo},
that $\hscomp$
is a free module over $H_{\Zt}^*(pt)$.
Then surjectivity of $\phi_2$ follows from surjectivity of $\phi$
using a formal argument identical to that of Corollary \ref{surjective}.
By Theorem \ref{tara} and Proposition \ref{formality},
$\hrs$ is a free module over $H_{\Zt}^*(pt)$,
therefore $\ker\phi_2$ is a free 
$H_{\Zt}^*(pt)$-module of rank $n$.
The fact that $u_i(x-u_i)\in \ker\phi_2$ follows from the
argument of Remark \ref{kernel}, hence we are done.
\end{proof}

The ring $\hscomp$ is therefore a deformation
of $\ost$ (over the base $\operatorname{Spec}\Zt[x]$)
that depends nontrivially on the affine structure of $\arr$,
rather than simply on the underlying matroid.

\begin{example}\label{first}
Consider the arrangements $\arr_a$ and $\arr_c$ in Figure 2(a) and 2(c).
By Theorem \ref{surj} and Example \ref{ts} we have
$$H^*_{\Zt}(\mathcal{M}(\arr_a);\Zt)\cong
\Zt [u_1,\ldots,u_4,x]\bigg/
\< \begin{array}{c}
u_1(x-u_1), u_2(x-u_2), u_3(x-u_3), u_4(x-u_4),\\
u_2u_3, u_1(x-u_2)u_4, u_1u_3u_4 
\end{array}\>$$
and
$$H^*_{\Zt}(\mathcal{M}(\arr_c);\Zt)\cong
\Zt [u_1,\ldots,u_4,x]\bigg/
\< \begin{array}{c}
u_1(x-u_1), u_2(x-u_2), u_3(x-u_3), u_4(x-u_4),\\
u_2u_3, (x-u_1)u_2(x-u_4), u_1u_3u_4
\end{array}\>.$$
The map $f:H^*_{\Zt}(\mathcal{M}(\arr_a);\Zt)\to H^*_{\Zt}(\mathcal{M}(\arr_b);\Zt)$
given by
$$f(u_1) = u_1+u_2, f(u_2) = u_2+u_3+x, f(u_3) = u_3,
f(u_4) = u_2+u_4,\text{ and }f(x)=x$$
is an isomorphism of graded $\Zt[x]$-algebras, hence the
ring $\hscomp$ is not a complete invariant of smooth, rational, 
affine arrangements up to combinatorial equivalence.\footnote{We thank
Graham Denham for finding this isomorphism.}
\end{example}

\begin{example}\label{second}
Now consider the arrangements $\arr_a'$ and $\arr_c'$ obtained from $\arr_a$ and $\arr_c$
by adding a vertical line on the far left, as shown below.
\begin{figure}[h]
\begin{center}
\psfrag{Ha'}{$\arr_a'$}
\psfrag{Hc'}{$\arr_c'$}
\includegraphics{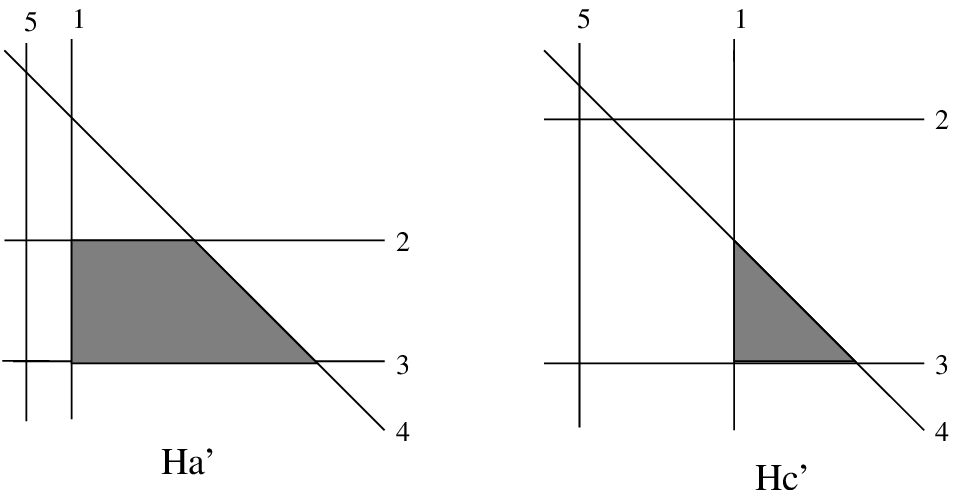}
\end{center}
\end{figure}
Again by Theorem \ref{surj}, we have
$$H^*_{\Zt}(\mathcal{M}(\arr_a');\Zt)\cong
\Zt [\vec{u},x]\bigg/
\< \begin{array}{c}
u_1(x-u_1), u_2(x-u_2), u_3(x-u_3), u_4(x-u_4),\\ u_5(x-u_5), 
u_2u_3, (x-u_1)u_5, u_1(x-u_2)u_4,\\ u_1u_3u_4, (x-u_2)u_4u_5, u_3u_4u_5 
\end{array}\>$$
and
$$H^*_{\Zt}(\mathcal{M}(\arr_c');\Zt)\cong
\Zt [\vec{u},x]\bigg/
\< \begin{array}{c}
u_1(x-u_1), u_2(x-u_2), u_3(x-u_3), u_4(x-u_4),\\ u_5(x-u_5),
u_2u_3, (x-u_1)u_5, (x-u_1)u_2(x-u_4),\\ u_1u_3u_4, (x-u_2)u_4u_5, u_3u_4u_5 
\end{array}\>.$$
We have used Macaulay 2 to check that the annihilator of the element 
$u_2\in H^*_{\Zt}(\mathcal{M}(\arr_a');\Zt)$ is generated by two linear elements
(namely $u_3$ and $x-u_2$) and nothing else, while there is no element
of $H^*_{\Zt}(\mathcal{M}(\arr_c');\Zt)$ with this property.
Hence the two rings are not isomorphic.
\end{example}
\end{section}

\footnotesize{

}

\end{spacing}
\end{document}